\documentclass[11pt]{article}
	
	\newcommand{\mainTitle}{Generators for vector spaces spanned by double zeta values with even weight}

	\newcommand{\emailAddressFst}{machide@nii.ac.jp}

\usepackage{graphicx}
\usepackage{geometry}               
\usepackage{ascmac}
\usepackage{amsmath}
\usepackage{amssymb}
\usepackage{amsthm}
\usepackage[hypertex, colorlinks=true, bookmarks=true]{hyperref}
\usepackage{arydshln}

	\DeclareMathOperator*{\OPlus}{\bigoplus}

	\newcommand{\nbk}[3]{#1#3#2}		
	\newcommand{\bgbk}[3]{\bigl{#1}#3\bigr{#2}}	
	\newcommand{\Bgbk}[3]{\Bigl{#1}#3\Bigr{#2}}			
	\newcommand{\bggbk}[3]{\biggl{#1}#3\biggr{#2}}			
	\newcommand{\Bggbk}[3]{\Biggl{#1}#3\Biggr{#2}}
	\newcommand{\autobk}[3]{\left#1#3\right#2}
	\newcommand{\nbkD}[5]{#1#2#5#3#4}		
	\newcommand{\bgbkD}[5]{\bigl{#1}\bigl{#2}#5\bigr{#3}\bigr{#4}}	
	\newcommand{\BgbkD}[5]{\Bigl{#1}\Bigl{#2}#5\Bigr{#3}\Bigr{#4}}	
	\newcommand{\bggbkD}[5]{\biggl{#1}\biggl{#2}#5\biggr{#3}\biggr{#4}}	
	\newcommand{\BggbkD}[5]{\Biggl{#1}\Biggl{#2}#5\Biggr{#3}\Biggr{#4}}	
	\newcommand{\autobkD}[5]{\left#1\left#2#5\right#3\right#4}	
	\newcommand{\mcbk}[4][?]{\ifx n#1\nbk{#2}{#3}{#4}\else\ifx b#1\bgbk{#2}{#3}{#4}\else\ifx B#1\Bgbk{#2}{#3}{#4}\else\ifx g#1\bggbk{#2}{#3}{#4}\else\ifx G#1\Bggbk{#2}{#3}{#4}\else\ifx a#1\autobk{#2}{#3}{#4}\else\ifx !#1{#4}\else#4\fi\fi\fi\fi\fi\fi\fi}
	\newcommand{\mcbkD}[4][?]{\ifx n#1\nbkD{#2}{#2}{#3}{#3}{#4}\else\ifx b#1\bgbkD{#2}{#2}{#3}{#3}{#4}\else\ifx B#1\BgbkD{#2}{#2}{#3}{#3}{#4}\else\ifx g#1\bggbkD{#2}{#2}{#3}{#3}{#4}\else\ifx G#1\BggbkD{#2}{#2}{#3}{#3}{#4}\else\ifx a#1\autobkD{#2}{#2}{#3}{#3}{#4}\else\ifx !#1{#4}\else#4\fi\fi\fi\fi\fi\fi\fi}
	\newcommand{\nsgsb}[1]{#1}		
	\newcommand{\bgsgsb}[1]{\big{#1}}	
	\newcommand{\Bgsgsb}[1]{\Big{#1}}			
	\newcommand{\bggsgsb}[1]{\bigg{#1}}			
	\newcommand{\Bggsgsb}[1]{\Bigg{#1}}
	\newcommand{\mcsgsb}[2][?]{\ifx n#1\nsgsb{#2}\else\ifx b#1\bgsgsb{#2}\else\ifx B#1\Bgsgsb{#2}\else\ifx g#1\bggsgsb{#2}\else\ifx G#1\Bggsgsb{#2}\else#2\fi\fi\fi\fi\fi}
	\newcommand{\myEqSpace}{\,} 	\newlength{\myEqSpaceLen} 	
	\newcommand{\mLt}[1]{\widetilde{#1}}

	\newcommand{\bkR}[2][a]{\mcbk[#1]{(}{)}{#2}}						
	\newcommand{\bkS}[2][a]{\mcbk[#1]{[}{]}{#2}}						
	\newcommand{\bkB}[2][a]{\mcbk[#1]{\{}{\}}{#2}}						
		
	\newcommand{\gam}{\gamma}

	\newcommand{\gpKleinF}[1][?]{V}

	\newcommand{\gpu}[1][?]{\ifx?#1e\else e_{#1}\fi}

	\newcommand{\vPack}[1][10]{\vspace{-#1pt}}
	
	\newcommand{\lnA}[1][]{&  &}
	\newcommand{\lnP}[1]{\myEqSpace#1\myEqSpace}

		\newcommand{\slnAH}[1][?]{\\}
		
	\newcommand{\pcstSpForRefThm}{\;}		\newcommand{\refHL}[2]{#1\pcstSpForRefThm\ref{#2}}		
	\newcommand{\refThm}[2][?]{\ifx?#1\refHL{Theorem}{#2}\else\ifx s#1\refHL{Theorems}{#2}\else{[argument error]}\fi\fi}
	\newcommand{\refProp}[2][?]{\ifx?#1\refHL{Proposition}{#2}\else\ifx s#1\refHL{Propositions}{#2}\else{[argument error]}\fi\fi}
	\newcommand{\refLem}[2][?]{\ifx?#1\refHL{Lemma}{#2}\else\ifx s#1\refHL{Lemmas}{#2}\else{[argument error]}\fi\fi}
	\newcommand{\refCor}[2][?]{\ifx?#1\refHL{Corollary}{#2}\else\ifx s#1\refHL{Corollaries}{#2}\else{[argument error]}\fi\fi}
	\newcommand{\refDef}[2][?]{\ifx?#1\refHL{Definition}{#2}\else\ifx s#1\refHL{Definitions}{#2}\else{[argument error]}\fi\fi}
	\newcommand{\refRem}[2][?]{\ifx?#1\refHL{Remark}{#2}\else\ifx s#1\refHL{Remarks}{#2}\else{[argument error]}\fi\fi}
	\newcommand{\refTab}[2][?]{\ifx?#1\refHL{Table}{#2}\else\ifx s#1\refHL{Tables}{#2}\else{[argument error]}\fi\fi}
	\newcommand{\refSec}[2][?]{\ifx?#1\refHL{Section}{#2}\else\ifx s#1\refHL{Sections}{#2}\else{[argument error]}\fi\fi}
		
	\newcommand{\Vc}[2][?]{\ifx ?#1\vec{#2}\else\ifx l#1\overrightarrow{#2}\else\ifx b#1{\bf#2}\else[error]\fi\fi\fi}
	
	\newcommand{\nopF}[3][?]{\ifx s#1#2/#3\else\ifx b#1(#2)/(#3)\else\ifx d#1\dfrac{#2}{#3}\else\ifx t#1\tfrac{#2}{#3}\else\frac{#2}{#3}\fi\fi\fi\fi}

	\newcommand{\myVopLetter}{$\cdot$} \newlength{\myVopLetterHeight}   

	\newcommand{\pw}[3][?]{\ifx!#3{#2}^{#3}\else#2^{#3}\fi}
	
	\newcommand{\pwR}[3][a]{\ifx!#1{\bkR[#1]{#2}}^{#3}\else\bkR[#1]{#2}^{#3}\fi}
	\newcommand{\pwB}[3][a]{\ifx!#1{\bkB[#1]{#2}}^{#3}\else\bkB[#1]{#2}^{#3}\fi}
	\newcommand{\pwS}[3][a]{\ifx!#1{\bkS[#1]{#2}}^{#3}\else\bkS[#1]{#2}^{#3}\fi}
					
	\newcommand{\mpIlett}{id}		\newcommand{\mpI}[1][?]{\ifx?#1\mpIlett\else \mpIlett_{#1}\fi}	
	
	\newcommand{\nSmO}[2][?]{\ifx l#1\sum\limits_{#2}\else\ifx t#1{\textstyle\sum\limits_{#2}}\else\sum_{#2}\fi\fi}
	\newcommand{\nSmT}[3][?]{\ifx l#1\sum\limits_{#2}^{#3}\else\if t#1{\textstyle\sum\limits_{#2}^{#3}}\else\sum_{#2}^{#3}\fi\fi}	
	\newcommand{\nSmN}[1][?]{\ifx l#1\sum\limits\else\ifx t#1{\textstyle\sum\limits}\else\sum\fi\fi}
	\newcommand{\pSm}[2][?]{\ifx t#1 \sum_{#2}^{\prime} \else \sideset{}{^\prime}\sum_{#2} \fi}
	\newcommand{\pSmT}[3][?]{\ifx t#1 \sum_{#2}^{\prime#3} \else \sideset{}{^\prime}\sum_{#2}^{#3} \fi}	
	\newcommand{\pSmN}[1][?]{\ifx t#1 \sum^{\prime} \else \sideset{}{^\prime}\sum \fi}
	\newcommand{\dSm}[2][?]{\ifx t#1 \sum_{#2}^{\dagger} \else \sideset{}{^\dagger}\sum_{#2} \fi}
	\newcommand{\dSmT}[3][?]{\ifx t#1 \sum_{#2}^{\dagger#3} \else \sideset{}{^\dagger}\sum_{#2}^{#3} \fi}	
	\newcommand{\dSmN}[1][?]{\ifx t#1 \sum^{\dagger} \else \sideset{}{^\dagger}\sum \fi}

	\newcommand{\nPd}[2][?]{\ifx l#1 \prod\limits_{#2} \else \prod_{#2} \fi}
	\newcommand{\nPdT}[3][?]{\ifx l#1 \prod\limits_{#2}^{#3} \else \prod_{#2}^{#3} \fi}

	\newcommand{\nOPs}[2][?]{\ifx l#1 \OPlus\limits_{#2} \else \OPlus_{#2} \fi}
	\newcommand{\nOPsT}[3][?]{\ifx l#1 \OPlus\limits_{#2}^{#3} \else \OPlus_{#2}^{#3} \fi}	
	\newcommand{\pOPs}[2][?]{\ifx t#1 \OPlus_{#2}^{\prime} \else \sideset{}{^\prime}\OPlus_{#2} \fi}
	\newcommand{\pOPsT}[3][?]{\ifx t#1 \OPlus_{#2}^{\prime#3} \else \sideset{}{^\prime}\OPlus_{#2}^{#3} \fi}

	\newcommand{\nIs}[2][?]{\ifx l#1 \bigcap\limits_{#2}\else\ifx b#1 \bigcap_{#2}\else{\textstyle\bigcap\limits_{#2}}\fi\fi}
	\newcommand{\nIsT}[3][?]{\ifx l#1 \bigcap\limits_{#2}^{#3}\else\ifx b#1 \bigcap_{#2}^{#3}\else{\textstyle\bigcap\limits_{#2}^{#3}}\fi\fi}	
	\newcommand{\pIs}[2][?]{\ifx t#1 \bigcap_{#2}^{\prime} \else \sideset{}{^\prime}\bigcap_{#2} \fi}
	\newcommand{\pIsT}[3][?]{\ifx t#1 \bigcap_{#2}^{\prime#3} \else \sideset{}{^\prime}\bigcap_{#2}^{#3} \fi}

	\newcommand{\nUn}[2][?]{\ifx l#1 \bigcup\limits_{#2}\else\ifx b#1 \bigcup_{#2}\else{\textstyle\bigcup\limits_{#2}}\fi\fi}
	\newcommand{\nUnT}[3][?]{\ifx l#1 \bigcup\limits_{#2}^{#3}\else\ifx b#1 \bigcup_{#2}^{#3}\else{\textstyle\bigcup\limits_{#2}^{#3}}\fi\fi}	
	\newcommand{\pUn}[2][?]{\ifx t#1 \bigcup_{#2}^{\prime} \else \sideset{}{^\prime}\bigcup_{#2} \fi}
	\newcommand{\pUnT}[3][?]{\ifx t#1 \bigcup_{#2}^{\prime#3} \else \sideset{}{^\prime}\bigcup_{#2}^{#3} \fi}

	\newcommand{\nLm}[2][?]{\ifx l#1 \lim\limits_{#2} \else \lim_{#2} \fi}
		
	\newcommand{\glcondEnvLineHead}[1]{ \ifx*#1 \begin{eqnarray*} \else \begin{eqnarray}  \label{#1} \fi }
	\newcommand{\glcondEnvLineTail}[1]{ \ifx*#1 \end{eqnarray*} \else \end{eqnarray} \fi }
	\newcommand{\glcondDis}[1]{\ifx d#1 \displaystyle \fi}
	
	\newcommand{\glcmdHLineCWiden}{\rule{0cm}{15pt}}	\newcommand{\glcdH}{\glcmdHLineCWiden}
	\newcommand{\lccondPar}[1]{\ifx#1p \\ \fi}

		\newcommand{\envMO}[2][*]{$\ifx d#1 \displaystyle \fi#2$}
		\newcommand{\envMT}[3][*]{$\ifx d#1 \displaystyle \fi#2=#3$}
		\newcommand{\envMTDef}[3][*]{$\ifx d#1 \displaystyle \fi#2:=#3$}
		\newcommand{\envMTPt}[4][*]{$\ifx d#1 \displaystyle \fi#3#2#4$}

		\newcommand{\envMTh}[4][*]{$\ifx d#1 \displaystyle \fi#2=#3=#4$}
		\newcommand{\envMThPt}[5][*]{$\ifx d#1 \displaystyle \fi#3#2#4#2#5$}
		\newcommand{\envMF}[5][*]{$\ifx d#1 \displaystyle \fi#2=#3=#4=#5$}
		\newcommand{\envMFPt}[6][*]{$\ifx d#1 \displaystyle \fi#3#2#4#2#5#2#6$}
		
	\newcommand{\envMLineT}[3][*]{ \ifx*#1 \begin{multline*} #2\lnP{=}#3\end{multline*} \else \begin{multline} \label{#1} #2\lnP{=}#3\end{multline} \fi }
	\newcommand{\envMLineTDef}[3][*]{ \ifx*#1 \begin{multline*} #2\lnP{:=}#3\end{multline*} \else \begin{multline} \label{#1} #2\lnP{:=}#3\end{multline} \fi }
	\newcommand{\envMLineTPt}[4][*]{ \ifx*#1 \begin{multline*} #3\lnP{#2}#4\end{multline*} \else \begin{multline} \label{#1} #3\lnP{#2}#4\end{multline} \fi }

		\newcommand{\envHLineT}[3][*]{ \glcondEnvLineHead{#1} #2&=&#3\glcondEnvLineTail{#1} }

	\newcommand{\matu}[1][?]{\ifx#1?I\else I_{#1}\fi}

	\newcommand{\ctG}[1][?]{\gam}

	\newcommand{\cTxT}[2]{\textcolor{#1}{#2}}
	
		\newcommand{\cTx}{\cTxT}
	\newcommand{\alTx}[1]{\cTx{red}{#1}}
	\newcommand{\rmTx}[1]{\cTx{blue}{#1}}
	\newcommand{\ntTx}[1]{\cTx{Magenta}{#1}}
	
	\newcommand{\cmoTx}[1]{\cTx{Gray}{#1}}
	\newcommand{\sTx}[2][?]{ \ifx t#1{\tiny #2} \else \ifx s#1{\scriptsize #2} \else \ifx f#1{\footnotesize #2} \else \ifx S#1{\small #2} \else \ifx n#1{\normalsize #2} \else \ifx l#1{\large #2} \else \ifx L#1{\Large #2} \else \ifx R#1{\LARGE #2} \else \ifx h#1{\huge #2} \else \ifx H#1{\Huge #2} \else \ifx ?#1 #2 \else #2 \fi\fi\fi\fi\fi\fi\fi\fi\fi\fi\fi }
	\newcommand{\bfTx}[1]{{\bf#1}}
	\newcommand{\raMTx}[3][?]{\raisebox{#2pt}[0pt][0pt]{$\ifx d#1\displaystyle\fi#3$}}

		\newcommand{\envHLineTCl}[3][a]{ \ifx a#1\alTx{\envHLineT{#2}{#3}} \else \ifx r#1 \rmTx{\envHLineT{#2}{#3}} \else\ifx n#1 \ntTx{\envHLineT{#2}{#3}}\else\ifx c#1 \cmoTx{\envHLineT{#2}{#3}} \else \text{[argument error]} \fi\fi\fi\fi \vPack[18] }

		\newcommand{\envHLineCSClPart}[8][?]{\ifx*#1 \begin{eqnarray*} \else \begin{eqnarray}  \label{#1}  \fi \alTx{#3}&#2&\alTx{#4}\\\glcdH#5&#2&#6\nonumber\\\glcdH\alTx{#7}&#2&\alTx{#8}\nonumber\glcondEnvLineTail{*}}
		
			\newcommand{\HLineCTCl}[3][?]{\alTx{#2}&=&\alTx{#3}\nonumber \ifx#1p \\\glcdH \fi}
			\newcommand{\HLineCTClDef}[3][?]{\alTx{#2}&:=&\alTx{#3}\nonumber \ifx#1p \\\glcdH \fi}
			\newcommand{\HLineCFCl}[5][?]{\alTx{#2}&=&\alTx{#3}\nonumber\\\glcdH#4&=&#5\nonumber \ifx#1p \\\glcdH \fi}
			\newcommand{\HLineCFClDef}[5][?]{\alTx{#2}&:=&\alTx{#3}\nonumber\\\glcdH#4&:=&#5\nonumber \ifx#1p \\\glcdH \fi}
			\newcommand{\HLineCSCl}[7][?]{\alTx{#2}&=&\alTx{#3}\nonumber\\\glcdH#4&=&#5\nonumber\\\glcdH\alTx{#6}&=&\alTx{#7}\nonumber\ifx#1p \\\glcdH \fi}
			\newcommand{\HLineCSClDef}[7][?]{\alTx{#2}&:=&\alTx{#3}\nonumber\\\glcdH#4&:=&#5\nonumber\\\glcdH\alTx{#6}&:=&\alTx{#7}\nonumber\ifx#1p \\\glcdH \fi}									
			\newcommand{\HLineCECl}[9][?]{\alTx{#2}&=&\alTx{#3}\nonumber\\\glcdH#4&=&#5\nonumber\\\glcdH\alTx{#6}&=&\alTx{#7}\nonumber\\\glcdH#8&=&#9\nonumber\ifx#1p \\\glcdH \fi}
			\newcommand{\HLineCEClDef}[9][?]{\alTx{#2}&:=&\alTx{#3}\nonumber\\\glcdH#4&:=&#5\nonumber\\\glcdH\alTx{#6}&:=&\alTx{#7}\nonumber\\\glcdH#8&:=&#9\nonumber\ifx#1p \\\glcdH \fi}
		
	\newcommand{\envCenter}[2][*]{\ifx*#1\begin{center}\else\begin{center}[#1]\fi #2\end{center}}
	\newcommand{\envFlushleft}[2][*]{\ifx*#1\begin{flushleft}\else\begin{flushleft}[#1]\fi #2\end{flushleft}}
	\newcommand{\envFlushright}[2][*]{\ifx*#1\begin{flushright}\else\begin{flushright}[#1]\fi #2\end{flushright}}
		
	\newcommand{\envItemIm}[2][*]{\ifx*#1\begin{itemize}\else\begin{itemize}[#1]\fi #2\end{itemize}}
	\newcommand{\envItemDp}[2][*]{\ifx*#1\begin{description}\else\begin{description}[#1]\fi #2\end{description}}
	\newcommand{\envItemEm}[2][*]{\ifx*#1\begin{enumerate}\else\begin{enumerate}[#1]\fi #2\end{enumerate}}

	\newcommand{\envMultCol}[3][*]{\ifx1#2#3\else\begin{multicols}{#2}\ifx*#1\else\mbox{}\vspace{-#1pt}\fi#3\end{multicols}\fi}
	
\theoremstyle{plain}
\newtheorem{theorem}{THEOREM}

\newtheorem{lemma}[theorem]{LEMMA}

\theoremstyle{definition}

\newtheorem{example}[theorem]{EXAMPLE}
\theoremstyle{remark}
\newtheorem{remark}[theorem]{REMARK}
\theoremstyle{plain}

\theoremstyle{definition}

\theoremstyle{remark}

\theoremstyle{plain}

\theoremstyle{definition}

\theoremstyle{remark}

\allowdisplaybreaks[4]

	\newcommand{\lccondBibitem}[3][]{ \if ?#2 \bibitem{#3} \else \bibitem[#2]{#3} \fi}
	\newcommand{\glcondEnvLineTailPd}[1]{.\ifx*#1 \end{eqnarray*} \else \end{eqnarray} \fi  }
	\newcommand{\glcondEnvLineTailCm}[1]{,\ifx*#1 \end{eqnarray*} \else \end{eqnarray} \fi }
	\newcommand{\prcondEnvEqSpHead}[1]{ \ifx*#1 \begin{equation*}[ERROR] \else \begin{equation}  \label{#1} \fi  }
	\newcommand{\prcondEnvEqSpTail}[1]{\ifx*#1 [ERROR]\end{equation*} \else \end{equation} \fi }

	\newcommand{\envProof}[2][?]{ \par\mbox{}\vspace{-5pt}\\ \ifx?#1\emph{Proof.}\else\emph{Proof of #1.}\fi \ #2 \hfill $\Box$\\ \par}
	
		\newcommand{\envHLineCFCmNme}[5][*]{\begin{eqnarray} #2&=&#3,\\\glcdH#4&=&#5 \glcondEnvLineTailCm{?} }
		\newcommand{\envHLineCFNmePd}[5][*]{\begin{eqnarray} #2&=&#3,\\\glcdH#4&=&#5 \glcondEnvLineTailPd{?} }
		\newcommand{\envHLineCFCmDefNme}[5][*]{\begin{eqnarray} #2&:=&#3,\\\glcdH#4&:=&#5 \glcondEnvLineTailCm{?} }
		\newcommand{\envHLineCFDefNmePd}[5][*]{\begin{eqnarray} #2&:=&#3,\\\glcdH#4&:=&#5 \glcondEnvLineTailPd{?} }

		\newcommand{\envHLineCFNmePdPt}[6][*]{\begin{eqnarray}#3&#2&#4,\\\glcdH#5&#2&#6\glcondEnvLineTailPd{?}}
		\newcommand{\envHLineCFCmNmePt}[6][*]{\begin{eqnarray}#3&#2&#4,\\\glcdH#5&#2&#6\glcondEnvLineTailCm{?}}
		\newcommand{\envHLineCFNmePdPte}[7][*]{\begin{eqnarray}#2&#3&#4,\\\glcdH#5&#6&#7\glcondEnvLineTailPd{?}}
		\newcommand{\envHLineCFCmNmePte}[7][*]{\begin{eqnarray}#2&#3&#4,\\\glcdH#5&#6&#7\glcondEnvLineTailCm{?}}

		\newcommand{\envHLineCSNmePd}[7][*]{\begin{eqnarray} #2&=&#3,\\\glcdH#4&=&#5,\\\glcdH#6&=&#7\glcondEnvLineTailPd{?}}
		\newcommand{\envHLineCSDefNmePd}[7][*]{\begin{eqnarray} #2&:=&#3,\\\glcdH#4&:=&#5,\\\glcdH#6&:=&#7\glcondEnvLineTailPd{?}}
		\newcommand{\envHLineCSCmNme}[7][*]{\begin{eqnarray} #2&=&#3,\\\glcdH#4&=&#5,\\\glcdH#6&=&#7\glcondEnvLineTailCm{?}}
		\newcommand{\envHLineCSCmDefNme}[7][*]{\begin{eqnarray} #2&:=&#3,\\\glcdH#4&:=&#5,\\\glcdH#6&:=&#7\glcondEnvLineTailCm{?}}

		\newcommand{\envHLineCSNmePdPt}[8][*]{\begin{eqnarray}#3&#2&#4,\\\glcdH#5&#2&#6,\\\glcdH#7&#2&#8\glcondEnvLineTailPd{?}}
		\newcommand{\envHLineCSCmNmePt}[8][*]{\begin{eqnarray}#3&#2&#4,\\\glcdH#5&#2&#6,\\\glcdH#7&#2&#8\glcondEnvLineTailCm{?}}
		\newcommand{\envHLineCSNmePdPte}[9][*]{\begin{eqnarray}#2&#3&#4,\\\glcdH#5&#6&#7,\\\glcdH#8&#2&#9\glcondEnvLineTailPd{?}}
		\newcommand{\envHLineCSCmNmePte}[9][*]{\begin{eqnarray}#2&#3&#4,\\\glcdH#5&#6&#7,\\\glcdH#8&#2&#9\glcondEnvLineTailCm{?}}
		
		\newcommand{\envHLineCENmePd}[9][*]{\begin{eqnarray} #2&=&#3,\\\glcdH#4&=&#5,\\\glcdH#6&=&#7,\\\glcdH#8&=&#9\glcondEnvLineTailPd{?}}
		\newcommand{\envHLineCEDefNmePd}[9][*]{\begin{eqnarray} #2&:=&#3,\\\glcdH#4&:=&#5,\\\glcdH#6&:=&#7,\\\glcdH#8&:=&#9\glcondEnvLineTailPd{?}}
		\newcommand{\envHLineCECmNme}[9][*]{\begin{eqnarray} #2&=&#3,\\\glcdH#4&=&#5,\\\glcdH#6&=&#7,\\\glcdH#8&=&#9\glcondEnvLineTailCm{?}}
		\newcommand{\envHLineCECmDefNme}[9][*]{\begin{eqnarray} #2&:=&#3,\\\glcdH#4&:=&#5,\\\glcdH#6&:=&#7,\\\glcdH#8&:=&#9\glcondEnvLineTailCm{?}}

			\newcommand{\pccondPaOMathEnvCmPdPar}[1]{\ifx#1p \\\glcdH \fi}
		\newcommand{\envMOCm}[2][*]{$\ifx d#1 \displaystyle \fi#2$,}
		\newcommand{\envMOPd}[2][*]{$\ifx d#1 \displaystyle \fi#2$.}
		
		\newcommand{\envMTCm}[3][*]{$\ifx d#1 \displaystyle \fi#2=#3$,}
		\newcommand{\envMTPd}[3][*]{$\ifx d#1 \displaystyle \fi#2=#3$.}
		\newcommand{\envMTCmDef}[3][*]{$\ifx d#1 \displaystyle \fi#2:=#3$,}
		\newcommand{\envMTDefPd}[3][*]{$\ifx d#1 \displaystyle \fi#2:=#3$.}
		\newcommand{\envMTCmPt}[4][*]{$\ifx d#1 \displaystyle \fi#3#2#4$,}
		\newcommand{\envMTPdPt}[4][*]{$\ifx d#1 \displaystyle \fi#3#2#4$.}

		\newcommand{\envMThCm}[4][*]{$\ifx d#1 \displaystyle \fi#2=#3=#4$,}
		\newcommand{\envMThPd}[4][*]{$\ifx d#1 \displaystyle \fi#2=#3=#4$.}
		\newcommand{\envMThCmPt}[5][*]{$\ifx d#1 \displaystyle \fi#3#2#4#2#5$,}
		\newcommand{\envMThPdPt}[5][*]{$\ifx d#1 \displaystyle \fi#3#2#4#2#5$.}
		\newcommand{\envMFCm}[5][*]{$\ifx d#1 \displaystyle \fi#2=#3=#4=#5$,}
		\newcommand{\envMFPd}[5][*]{$\ifx d#1 \displaystyle \fi#2=#3=#4=#5$.}
		\newcommand{\envMFCmPt}[6][*]{$\ifx d#1 \displaystyle \fi#3#2#4#2#5#2#6$,}
		\newcommand{\envMFPdPt}[6][*]{$\ifx d#1 \displaystyle \fi#3#2#4#2#5#2#6$.}

		\newcommand{\envHLineCFCmNm}[5][*]{ \begin{equation}\begin{split} \ifx*#1 \text{[ERROR;need label name]} \else \label{#1} \fi #2&\lnP{=}#3,\\#4&\lnP{=}#5, \end{split}\end{equation} }
		\newcommand{\envHLineCFNm}[5][*]{ \begin{equation}\begin{split} \ifx*#1 \text{[ERROR;need label name]} \else \label{#1} \fi #2&\lnP{=}#3\\#4&\lnP{=}#5, \end{split}\end{equation} }
		\newcommand{\envHLineCFNmPd}[5][*]{ \begin{equation}\begin{split} \ifx*#1 \text{[ERROR;need label name]} \else \label{#1} \fi #2&\lnP{=}#3,\\#4&\lnP{=}#5. \end{split}\end{equation} }
		\newcommand{\envHLineCFCmDefNm}[5][*]{ \begin{equation}\begin{split} \ifx*#1 \text{[ERROR;need label name]} \else \label{#1} \fi #2&\lnP{:=}#3,\\#4&\lnP{:=}#5, \end{split}\end{equation} }
		\newcommand{\envHLineCFDefNm}[5][*]{ \begin{equation}\begin{split} \ifx*#1 \text{[ERROR;need label name]} \else \label{#1} \fi #2&\lnP{:=}#3\\#4&\lnP{:=}#5, \end{split}\end{equation} }
		\newcommand{\envHLineCFDefNmPd}[5][*]{ \begin{equation}\begin{split} \ifx*#1 \text{[ERROR;need label name]} \else \label{#1} \fi #2&\lnP{:=}#3,\\#4&\lnP{:=}#5. \end{split}\end{equation} }
		\newcommand{\envHLineCSCmNm}[7][*]{ \begin{equation}\begin{split} \ifx*#1 \text{[ERROR;need label name]} \else \label{#1} \fi #2&\lnP{=}#3,\\#4&\lnP{=}#5,\\#6&\lnP{=}#7 \end{split}\end{equation} }
		\newcommand{\envHLineCSNm}[7][*]{ \begin{equation}\begin{split} \ifx*#1 \text{[ERROR;need label name]} \else \label{#1} \fi #2&\lnP{=}#3\\#4&\lnP{=}#5\\#6&\lnP{=}#7 \end{split}\end{equation} }
		\newcommand{\envHLineCSNmPd}[7][*]{ \begin{equation}\begin{split} \ifx*#1 \text{[ERROR;need label name]} \else \label{#1} \fi #2&\lnP{=}#3,\\#4&\lnP{=}#5,\\#6&\lnP{=}#7. \end{split}\end{equation} }
		\newcommand{\envHLineCSCmDefNm}[7][*]{ \begin{equation}\begin{split} \ifx*#1 \text{[ERROR;need label name]} \else \label{#1} \fi #2&\lnP{:=}#3,\\#4&\lnP{:=}#5,\\#6&\lnP{:=}#7 \end{split}\end{equation} }
		\newcommand{\envHLineCSDefNm}[7][*]{ \begin{equation}\begin{split} \ifx*#1 \text{[ERROR;need label name]} \else \label{#1} \fi #2&\lnP{:=}#3\\#4&\lnP{:=}#5\\#6&\lnP{:=}#7 \end{split}\end{equation} }
		\newcommand{\envHLineCSDefNmPd}[7][*]{ \begin{equation}\begin{split} \ifx*#1 \text{[ERROR;need label name]} \else \label{#1} \fi #2&\lnP{:=}#3,\\#4&\lnP{:=}#5,\\#6&\lnP{:=}#7. \end{split}\end{equation} }
		\newcommand{\envHLineCECmNm}[9][*]{ \begin{equation}\begin{split} \ifx*#1 \text{[ERROR;need label name]} \else \label{#1} \fi #2&\lnP{=}#3,\\#4&\lnP{=}#5,\\#6&\lnP{=}#7,\\#8&\lnP{=}#9,  \end{split}\end{equation} }
		\newcommand{\envHLineCENm}[9][*]{ \begin{equation}\begin{split} \ifx*#1 \text{[ERROR;need label name]} \else \label{#1} \fi #2&\lnP{=}#3\\#4&\lnP{=}#5\\#6&\lnP{=}#7\\#8&\lnP{=}#9  \end{split}\end{equation} }
		\newcommand{\envHLineCENmPd}[9][*]{ \begin{equation}\begin{split} \ifx*#1 \text{[ERROR;need label name]} \else \label{#1} \fi #2&\lnP{=}#3,\\#4&\lnP{=}#5,\\#6&\lnP{=}#7,\\#8&\lnP{=}#9.  \end{split}\end{equation} }
		\newcommand{\envHLineCECmDefNm}[9][*]{ \begin{equation}\begin{split} \ifx*#1 \text{[ERROR;need label name]} \else \label{#1} \fi #2&\lnP{:=}#3,\\#4&\lnP{:=}#5,\\#6&\lnP{:=}#7,\\#8&\lnP{:=}#9,  \end{split}\end{equation} }
		\newcommand{\envHLineCEDefNm}[9][*]{ \begin{equation}\begin{split} \ifx*#1 \text{[ERROR;need label name]} \else \label{#1} \fi #2&\lnP{:=}#3\\#4&\lnP{:=}#5\\#6&\lnP{:=}#7\\#8&\lnP{:=}#9  \end{split}\end{equation} }
		\newcommand{\envHLineCEDefNmPd}[9][*]{ \begin{equation}\begin{split} \ifx*#1 \text{[ERROR;need label name]} \else \label{#1} \fi #2&\lnP{:=}#3,\\#4&\lnP{:=}#5,\\#6&\lnP{:=}#7,\\#8&\lnP{:=}#9.  \end{split}\end{equation} }
	\newcommand{\envMLineTPd}[3][*]{ \ifx*#1 \begin{multline*} #2\lnP{=}#3.\end{multline*} \else \begin{multline} \label{#1} #2\lnP{=}#3.\end{multline} \fi }
	\newcommand{\envMLineTCm}[3][*]{ \ifx*#1 \begin{multline*} #2\lnP{=}#3,\end{multline*} \else \begin{multline} \label{#1} #2\lnP{=}#3,\end{multline} \fi }
	\newcommand{\envMLineTDefPd}[3][*]{ \ifx*#1 \begin{multline*} #2\lnP{:=}#3.\end{multline*} \else \begin{multline} \label{#1} #2\lnP{:=}#3.\end{multline} \fi }
	\newcommand{\envMLineTCmDef}[3][*]{ \ifx*#1 \begin{multline*} #2\lnP{:=}#3,\end{multline*} \else \begin{multline} \label{#1} #2\lnP{:=}#3,\end{multline} \fi }

	\newcommand{\envItemEmPa}[2][*]{\vspace{-5pt}\ifx*#1\begin{enumerate}\else\begin{enumerate}[#1]\fi \renewcommand{\itemsep}{0pt}\renewcommand{\itemindent}{0pt}\renewcommand{\theenumi}{\bfTx{(\roman{enumi})}}#2\end{enumerate}}

	\DeclareFontEncoding{OT2}{}{}
	\DeclareFontSubstitution{OT2}{cmr}{m}{n}
	\DeclareFontFamily{OT2}{cmr}{\hyphenchar\font45}
	\DeclareFontShape{OT2}{cmr}{m}{n}{<5><6><7><8><9>gen*wncyr <10><10.95><12><14.4><17.28><20.74><24.88>wncyr10}{}
	\DeclareFontShape{OT2}{cmr}{b}{n}{<5><6><7><8><9>gen*wncyb<10><10.95><12><14.4><17.28><20.74><24.88>wncyb10}{}
	\DeclareMathAlphabet{\mathcyr}{OT2}{cmr}{m}{n}
	\DeclareMathAlphabet{\mathcyb}{OT2}{cmr}{b}{n}
	\SetMathAlphabet{\mathcyr}{bold}{OT2}{cmr}{b}{n}

	\newcommand{\sh}{\mathcyr{sh}}

	\newcommand{\fcZBlettHelp}{\bullet}
	\newcommand{\fcZBN}[1][?]{\zeta^\fcZBlettHelp}			
	\newcommand{\fcZHN}[1][?]{\zeta^*}									
	\newcommand{\fcZRN}[1][?]{\zeta^\ddagger}								
	\newcommand{\fcZSN}[1][?]{\zeta^\sh}							
								
	\newcommand{\LetHelpFcCh}{\chi}
	\newcommand{\fcChBN}[1][?]{\LetHelpFcCh^\fcZBlettHelp}			
	\newcommand{\fcChHN}[1][?]{\LetHelpFcCh^*}									
	\newcommand{\fcChRN}[1][?]{\LetHelpFcCh^\ddagger}								
	\newcommand{\fcChSN}[1][?]{\LetHelpFcCh^\sh}							

	\newcommand{\fcChBmN}[1][?]{\mLt{\LetHelpFcCh}^\fcZBlettHelp}			
	\newcommand{\fcChHmN}[1][?]{\mLt{\LetHelpFcCh}^*}							
	\newcommand{\fcChRmN}[1][?]{\mLt{\LetHelpFcCh}^\ddagger}								
	\newcommand{\fcChSmN}[1][?]{\mLt{\LetHelpFcCh}^\sh}							
								
			\newcommand{\ggfcZlett}{\mathfrak{\bLt{Z}}}	
				
												\newcommand{\ggfcZbb}[1][?]{\ggfcZlett}
										\newcommand{\ggfcZHbb}[1][?]{\ggfcZlett^*}
							\newcommand{\ggfcZRbb}[1][?]{\ggfcZlett^\ddagger}
									\newcommand{\ggfcZSbb}[1][?]{\ggfcZlett^\sh}

\newcommand{\lgpeD}[2][?]{\ifx#1?\mathsf{D}_{#2}\else\mathsf{D}_{#2}^{(#1)}\fi}
\newcommand{\lgpeR}[1][?]{\ifx#1?\mathsf{R}\else\mathsf{R}^{(#1)}\fi}
\newcommand{\lgpeRO}[2][?]{\ifx#1?\mathsf{R}_{#2}\else\mathsf{R}_{#2}^{(#1)}\fi}
\newcommand{\lgpeRm}[1][?]{\ifx#1?\mLt{\mathsf{R}}\else\mLt{\mathsf{R}}^{(#1)}\fi}
\newcommand{\lgpeP}[1][?]{\ifx#1?\mathsf{P}\else\mathsf{P}^{(#1)}\fi}
\newcommand{\lgpePO}[2][?]{\ifx#1?\mathsf{P}_{#2}\else\mathsf{P}_{#2}^{(#1)}\fi}
\newcommand{\lgpePi}[1][?]{\ifx#1?\mathsf{P}^{-1}\else(\mathsf{P}^{(#1)})^{-1}\fi}
\newcommand{\lgpePOi}[2][?]{\ifx#1?(\mathsf{P}_{#2})^{-1}\else(\mathsf{P}_{#2}^{(#1)})^{-1}\fi}

	\newcommand{\lsstSh}[1][?]{\ifx#1?{\mathsf{S}}_{\sh}\else{\mathsf{S}}_{\sh}^{(#1)}\fi}
	\newcommand{\lsstShT}[2][?]{\ifx#1?{\mathsf{S}}_{\sh,#2}\else{\mathsf{S}}_{\sh,#2}^{(#1)}\fi}

	\renewcommand{\matu}[1][?]{\ifx#1?\mathsf{I}\else {\mathsf{I}^{(#1)}}\fi}
	\renewcommand{\gpu}[2][?]{\ifx ?#1e^{(#2)}\else e\fi}
	\newcommand{\pbfL}[1][?]{\bfTx{l}}
	\newcommand{\pbfK}[1][?]{\bfTx{k}}
	\newcommand{\pbfX}[1][?]{\bfTx{x}}

			\renewcommand{\ggfcZlett}{F}

											\renewcommand{\ggfcZbb}[1][?]{\ggfcZlett}
										\renewcommand{\ggfcZHbb}[1][?]{\ggfcZlett^*}
						\renewcommand{\ggfcZRbb}[1][?]{\ggfcZlett^\ddagger}
									\renewcommand{\ggfcZSbb}[1][?]{\ggfcZlett^\sh}

\allowdisplaybreaks[4]
	\title{\mainTitle}
	\date{}

\begin{document}
\maketitle
\renewcommand{\thefootnote}{\fnsymbol{footnote}}
\footnote[0]{e-mail : \emailAddressFst, machide.t@gmail.com}
\renewcommand{\thefootnote}{\arabic{footnote}}\setcounter{footnote}{0}
\vPack[30]

\begin{abstract}
Let $\mathcal{DZ}_k$ be the $\mathbb{Q}$-vector space spanned by
	double zeta values with weight $k$,
	and $\mathcal{DM}_k$ be its quotient space divided by 
	the space $\mathcal{PZ}_k$ spanned by 
	the zeta value $\zeta(k)$ and 
	products of two zeta values with total weight $k$.
When $k$ is even, an upper bound 
	for the dimension of $\mathcal{DM}_k$ is known.
By adding the dimensions of $\mathcal{DM}_k$ and $\mathcal{PZ}_k$,
	an upper bound of $\mathcal{DZ}_k$
	which equals $k/2$ minus 
	the dimension of the space of modular forms of weight $k$ on the modular group 
	is given.
In this note,
	we obtain some specific sets of  generators for $\mathcal{DM}_k$
	which represent the upper bound.
These yield the corresponding sets and the upper bound for $\mathcal{DZ}_k$.
\end{abstract}


\section{Introduction and main theorem}
In recent years, the multiple zeta values 
have appeared in various contexts in mathematics and physics (cf. \cite{AKNO,BBBL,zagier2}).
The double zeta values, which are the multiple zeta values with depth $2$
and also called Euler sums or Euler-Zagier sums, are defined by
\begin{equation*}
\zeta (q,p) :=
\sum_{n > m> 0}
\frac{1}{n^{q} m^{p}}
\end{equation*}
for integers $q \geq 2$ and $p \geq 1$.
The integer $k = p+q$ is the weight of $\zeta (q,p)$ by definition.
These values
 go back to L. Euler \cite{euler1},
and were revisited by N. Nielsen \cite{nielsen1}, L. Tornheim \cite{tornheim1}
and L.J. Mordell \cite{mordell1}.
Euler discovered that
 the double zeta value $\zeta (q,p)$ is 
 a $\mathbb{Q}$-linear combination of 
 $\zeta(k)$ and  $\zeta(j) \zeta(k-j) \ (2 \leq j \leq k-2)$  
when its weight $k$ is odd
(see \cite[Introduction]{BBG} and \cite[Theorem 1]{HWZ} for an explicit combination).
Furthermore, he gave the formula 
\begin{equation}\label{1.Euler_formula}
2\zeta (k-1,1)=
(k-1) \zeta (k) - \sum_{j=2}^{k-2} \zeta(j) \zeta(k-j)
\end{equation}
for any integer $k \geq 3$.
Here the real numbers $\zeta(j) := \sum_{m =1}^{\infty} 1/m^j$ are the usual zeta values.

Let $\mathcal{DZ}_k $ 
be the $\mathbb{Q}$-vector space spanned by double zeta values with weight $k$,
and $\mathcal{DM}_k$ be its quotient space
$\mathcal{DZ}_k/ \mathcal{PZ}_k $, 
where $\mathcal{PZ}_k$ denotes
the space spanned by
the zeta value $\zeta (k)$ and 
products 
$\{ \zeta(j) \zeta(k-j)  \vert \ 2 \leq j \leq k-2\}$ of two zeta values with total weight $k$.
When $k$ is odd,
 the combinations discovered by
 Euler give generators $\{\zeta(k)\} \cup \{ \zeta(j) \zeta(k-j)  \vert \ 2 \leq j \leq (k-1)/2\}$
  for the space $\mathcal{DZ}_k $.
The generators imply that 
the integer $(k-1)/2$ is an upper bound for its dimension,
and $\mathcal{DM}_k$
is the null space.
Thus the spaces $\mathcal{DZ}_k $ and $\mathcal{DM}_k$ for $k$ odd are simpler
than those for $k$ even as we see below.

Let  $[x]$ be the greatest integer not exceeding $x$,
and $M_k$ be the space of 
modular forms of weight $k$ 
on the modular group $\mathrm{SL}(2, \mathbb{Z})$.
When $k$ is even,
D. Zagier has obtained an upper bound $[(k-2)/6]$
for the dimension of 
the space $\mathcal{DM}_k$ (see \cite[Section 8 and Appendix]{IKZ}, \cite{zagier1}, \cite[Section 8]{zagier2}).
Since the $[(k+2)/4]$ numbers $\{\zeta(k)\} \cup \{\zeta(j) \zeta(k-j) \vert \ 2 \leq j \leq k/2,\ \text{$j$ $odd$}\}$
generate the space $\mathcal{PZ}_k$
by Euler's result $\zeta (2j) \in \mathbb{Q} \pi^{2j}$,
the integer $[(k+2)/4] + [(k-2)/6] = k/2 - \mathrm{dim}\  M_k$
is an upper bound for the dimension of the space $\mathcal{DZ}_k$.
This upper bound
has been also obtained in the paper \cite[Theorem 2 and 3]{GKZ} of H. Gangl \emph{et al.}
by using surprising connections between the structure of $\mathcal{DZ}_k $ and 
that of $M_k$.
It should be noted that there are related works 
\cite{ihara, goncharov, schneps}
(see Remark in \cite[Introduction]{GKZ}).
However they have not discovered specific generators 
for $\mathcal{DM}_k$ and $\mathcal{DZ}_k$
which give the upper bounds.

The purpose of this note
is to prove the following theorem, more precisely to give $2^{[(k-2)/6]}$ specific sets of $[(k-2)/6]$ generators 
for the space $\mathcal{DM}_k$ when $k$ is even.
These sets clearly yield at least $2^{[(k-2)/6]}$ specific sets of 
$(k/2 - \mathrm{dim}\  M_k)$ generators for the space $\mathcal{DZ}_k$.
\begin{theorem}
\label{1.Theorem1}
Let $k$ be an even integer with $k \geq 2$.
For any choice of $\varepsilon_1, \ldots, \varepsilon_{[(k-2)/6]} \in \{ 0, 1 \}$,
the $[(k-2)/6]$ double zeta values
$\{\zeta(2j+\varepsilon_j ,k- 2j -\varepsilon_j) \vert \ 1 \leq j \leq [(k-2)/6] \}$
generate  
the $\mathbb{Q}$-vector space $\mathcal{DM}_k$.
The integer $[(k-2)/6]$ is consequently an upper bound for its dimension.
\end{theorem}


Since $[(k-2)/6]= [[(k-1)/3]/2]$ if $k$ is even,
Theorem \ref{1.Theorem1} with $\varepsilon_1= \cdots = \varepsilon_{[(k-2)/6]} =0$ says that 
the space $\mathcal{DM}_k$ is generated by
\begin{equation*}\label{1.a}
\{\zeta(2j ,k- 2j) \vert \ 1 \leq j \leq [\frac{k-2}{6}] \}
=
\{\zeta(j,k-j) \vert \ 2 \leq j \leq [\frac{k-1}{3}],\ \text{$j$ $even$} \}.
\end{equation*}

Thus $\mathcal{DM}_k$ are spanned by $\zeta(even, even)$'s
with weight $k$.

We know from \cite[Theorem 2, 3]{GKZ} that
 $\zeta(odd, odd)$'s with weight $k$ generate
 $\mathcal{DZ}_k$, and furthermore 
satisfy at least $(\mathrm{dim}\ M_k -1)$
linearly independent relations.
As an application of Theorem \ref{1.Theorem1},
we give a specific set of $(k/2 - \mathrm{dim}\  M_k)$ generators written in terms of $\zeta(odd, odd)$'s
for $\mathcal{DZ}_k$,
which reproduces the known results.
\begin{theorem}\label{1.Theorem2}
$\mathrm{(}$cf. \cite[Theorem 2, 3]{GKZ}$\mathrm{)}$
Let $k$ be an even integer with $k \geq 2$.\\
{\bf (i)}
 The following $(k/2 - \mathrm{dim}\  M_k)$ linear combinations of $\zeta(odd, odd)$'s
generate the $\mathbb{Q}$-vector space $\mathcal{DZ}_k$.
\begin{multline*}
\{\zeta(j,k-j) \vert \ 3 \leq j \leq [\frac{k+2}{3}],\ \text{j odd } \}
 \cup
\{\zeta(k-j,j) \vert \ 1 \leq j \leq [\frac{k+2}{3}],\ \text{j odd } \} \\
 \cup
\{\zeta(j,k-j)+ \zeta(k-j,j) \vert \ 
[\frac{k+5}{3}] \leq j \leq \frac{k}{2},\ \text{j odd } \}.
\end{multline*}
{\bf (ii)}
For any odd integer $i$ with $[(k+5)/3] \leq i  < k/2$,
there is a $\mathbb{Q}$-linear relation among $\zeta(odd, odd)$'s such that
\begin{equation*}
\sum_{j =3 \atop j\ odd}^{k-1} 
c_j \zeta (j, k-j)
 =
0
\qquad
(c_j \in \mathbb{Q})
\end{equation*}
where
$c_i \neq c_{k-i}$ 
and
$c_j = c_{k-j}$
if $j$ is an odd integer with $[(k+5)/3] \leq j  < k/2$ and $j \neq i$.
In particular, the number 
of linearly independent relations of the above kind equals $(\mathrm{dim}\ M_k -1)$.
\end{theorem}

\begin{remark}
The equations obtained in Theorem \ref{1.Theorem2}(ii) are presumably unique up to scaling,
	but the uniqueness can not be proved;
We do not know even whether all double zeta values $\zeta (j,k-j)$ are not rational.
\end{remark}

In the proof of Theorem \ref{1.Theorem1} below,
we shall use only relations among double zeta values and Tornheim double series 
which are defined by
\begin{equation}
T(r, q,p) 
:=
\sum_{n,m> 0}
\frac{1}{(n+m)^r n^q m^p}
\end{equation}
for integers $r,q,p \geq 1$.
We have reversed the order of arguments in $T(r, q,p)$ compared to the standard definition
in order to keep 
the relationship 
$T (r, 0, p) = \zeta (r, p)$
 with the summation convention 
used for the double zeta values.

Some basic properties of $T(r, q,p)$ are the following. 
\begin{align*}
T(r,q,p) &= T(r,p,q), & \\
T(r,q,p) &= T(r+1,q-1,p) + T(r+1, q, p-1) &( p,q \geq 1), \\
T(r, 0, p ) &= \zeta (r, p)  &(r \geq 2). 
\end{align*}
In the next and final section, 
	we prove Theorem \ref{1.Theorem1} and \ref{1.Theorem2},
	and 
	give some examples of equations obtained in Theorem \ref{1.Theorem2}(ii).

\section{Proofs of Theorem \ref{1.Theorem1} and \ref{1.Theorem2}}
Throughout this section, we suppose that
$k \geq 2$ is even.
Let $\mathcal{DZ}_k(r)$ be 
the $\mathbb{Q}$-vector space generated by  
 double zeta values $\{\zeta(j,k-j) \vert \ 2 \leq j \leq r\}$
for any integer $r$ with $1 \leq r \leq k-1$.
The space $\mathcal{DZ}_k(r)$ stands for the null space if $r=1$,
and equals $\mathcal{DZ}_k$ if $r=k-1$.
To give a proof of Theorem \ref{1.Theorem1}, we need some relations among double zeta values and Tornheim double series.
\begin{lemma} \label{2_Lem}
Let $p, q, r$ be integers with $k = p+ q+ r$.\\
{\bf (i)} If $p, q, r \geq 1$, then
\begin{equation}\label{2.Lemma1i}
(-1)^rT(r,q,p) + (-1)^qT(q,p,r) + (-1)^p T(p,r,q) \in \mathbb{Q} \zeta(k) \subset \mathcal{PZ}_k.
\end{equation}
{\bf (ii)} If $p, q, r \geq 1$, then
\begin{equation}\label{2.Lemma1iia}
T(r,q,p) \in \mathcal{DZ}_k(r) + \mathcal{PZ}_k.
\end{equation}
In particular, under the extra assumption that $r \geq 2$,
we have
\begin{equation}\label{2.Lemma1iib}
T(r,q,p)- (-1)^p \zeta(r, k-r) 
\in
\mathcal{DZ}_k(r-1) +  \mathcal{PZ}_k.
\end{equation}
{\bf (iii)} If $3 \leq r \leq k-3$ and $r$ is odd, then
\begin{equation}\label{2.Lemma1iii}
2\zeta(r, k-r) + (k-r) \zeta(r-1, k-r+1)
\in
\mathcal{DZ}_k(r-2) +  \mathcal{PZ}_k.
\end{equation}
\end{lemma}
\emph{Proof.}\ \ 
From \cite{EM, nakamura1, tsumura1, zagier3} (see also \cite{MNOT, mordell1, SS}), 
we obtain (\ref{2.Lemma1i}).
If $p \geq 2$, $q \geq 0$ and $r \geq 1$ with $q +r \geq 2$,
then, by \cite[Corollary 2.4]{boyadzhiev1} and (\ref{1.Euler_formula}),
\begin{align*}
T(r, q, p) 
-
(-1)^p
\sum_{j =1}^{r-1} \binom{p+r-j-2}{p-1} \zeta (j+1, k-j-1) 
\in
\mathcal{PZ}_k.
\end{align*}
Note that
$S(p,q)$ and $\zeta ( 1 )$ in \cite{boyadzhiev1} equal $\zeta (q,p) + \zeta (p+q)$ and $0$ respectively.
By using the above formula,
we deduce (\ref{2.Lemma1iia}) and (\ref{2.Lemma1iib})
from
$T(r, k-r-1,1) = T(r-1,k-r-1,2) - T(r, k-r-2,2)$,
and (\ref{2.Lemma1iii})
from
$T(r,0,k-r) = \zeta (r, k-r)$.
\qed \\ \\
We prove Theorem \ref{1.Theorem1}.\\ \\
\emph{Proof of Theorem \ref{1.Theorem1}.}\ \ 
Firstly we show that $\mathcal{DZ}_k = \mathcal{DZ}_k([(k-1)/3])+ \mathcal{PZ}_k$.
Let $r$ be an integer with $[k/3] +1 \leq r \leq k-2$.
If $k \not\equiv 2 \pmod{3}$ or $r \geq [k/3] +2$, then
there are integers $1\leq p, q \leq [k/3]$ such that $p+q+r=k$.
We see from (\ref{2.Lemma1i}) and (\ref{2.Lemma1iia}) that
$T(r,q,p) \in \mathcal{DZ}_k([k/3])+ \mathcal{PZ}_k$.
If $k \equiv 2 \pmod{3}$ and $r = [k/3] +1$,
then $p+q+r=k$ by setting $p= [k/3] $ and $q =[k/3] +1$.
Thus we also get by (\ref{2.Lemma1i}) and (\ref{2.Lemma1iia}) 
that $T(r,q,p) \in \mathcal{DZ}_k([k/3])+ \mathcal{PZ}_k$.
These together with (\ref{2.Lemma1iib}) imply
that, for any integer $r$ with $[k/3] +1 \leq r \leq k-2$,
we have $\zeta (r, k-r) \in \mathcal{DZ}_k(r-1)+ \mathcal{PZ}_k$,
or $\mathcal{DZ}_k(r) \subset \mathcal{DZ}_k(r-1)+ \mathcal{PZ}_k$
which is also true when $r=k-1$ by (\ref{1.Euler_formula}).
By using induction on $r$, starting at $r = [k/3] + 1$, we can obtain 
$\mathcal{DZ}_k  \subset \mathcal{DZ}_k([k/3])+ \mathcal{PZ}_k$.
If $k \equiv 0 \pmod{3}$, we see from (\ref{2.Lemma1i})
that $T(k/3,k/3,k/3) \in \mathcal{PZ}_k$
(which was obtained by Mordell \cite{mordell1}), 
and from (\ref{2.Lemma1iib})
that $\zeta(k/3, 2k/3) \in \mathcal{DZ}_k((k/3) -1)+ \mathcal{PZ}_k$.
Thus we conclude that
$\mathcal{DZ}_k = \mathcal{DZ}_k([(k-1)/3])+ \mathcal{PZ}_k$.

Let $\mathcal{X}_{k; \vec{\varepsilon}}(r)$
denote the space spanned by 
$\{\zeta(2j+\varepsilon_j ,k- 2j -\varepsilon_j) \vert \ 1 \leq j \leq r \}$
for any integer $r$ with $1 \leq r \leq [(k-2)/6]$.
We see from (\ref{2.Lemma1iii})
that $\mathcal{DZ}_k(2r+1) \subset 
\mathcal{X}_{k; \vec{\varepsilon}}(r)
+\mathcal{DZ}_k(2r-1) + \mathcal{PZ}_k$.
By induction on $r$, starting at $r = 1$,
we get 
$\mathcal{DZ}_k(2[(k-2)/6]+1) \subset 
\mathcal{X}_{k; \vec{\varepsilon}}([(k-2)/6]) + \mathcal{PZ}_k$.
Since $\mathcal{DZ}_k([(k-1)/3]) \subset \mathcal{DZ}_k(2[(k-2)/6]+1)$,
we obtain $\mathcal{DZ}_k = 
\mathcal{X}_{k; \vec{\varepsilon}}([(k-2)/6]) + \mathcal{PZ}_k$.
This completes the proof.
\qed \\ \\
Finally we give a proof of Theorem \ref{1.Theorem2}.\\ \\
\emph{Proof of Theorem \ref{1.Theorem2}.}\ \ 
From Theorem \ref{1.Theorem1} 
with $\varepsilon_1= \ldots= \varepsilon_{[(k-2)/6]} =1$
and the harmonic relation 
$\zeta(p) \zeta(q)
=
\zeta(p,q) +
\zeta(q,p) + \zeta (p+q)$,
we see that
the $(k/2 - \mathrm{dim}\  M_k)$ real numbers 
$\{\zeta(j,k-j) \vert \ 3 \leq j \leq [(k+2)/3],\ \text{$j$ $odd$ } \}
\cup
\{\zeta(j,k-j)+ \zeta(k-j,j) \vert \ 
2 \leq j \leq k/2,\ \text{$j$ $odd$ } \}
\cup
\{ \zeta (k) \}
$
generate the space $\mathcal{DZ}_k$.
This together with the formula \cite[Theorem 1]{GKZ}
\begin{equation*}
\sum_{j =3 \atop j\ odd}^{k-1} 
\zeta (j, k-j)
=
\frac{1}{4} \zeta (k)
\end{equation*}
yields (i).
We next prove (ii).
By (i),
the fact that
$
\zeta(i, k-i) 
\in
\mathcal{DZ}_k$
for any odd integer $i$ with $[(k+5)/3] \leq i  < k/2$
implies
the nontrivial equations stated in (ii).
If these are linearly dependent, 
	there are rational numbers 
	$d_i\in\mathbb{Q}$
	and
	$\mathbb{Q}$-linear combinations 
	$
	A^{(i)}
	=
	\sum\limits_{j =3 \atop j\ odd}^{k-1}  c_j^{(i)} \zeta (j, k-j)
	$
	for odd integers $i$ with $[(k+5)/3] \leq i  < k/2$
	such that
	$d_h \neq 0$ for some odd integer $h$ 
	and
\begin{equation} \label{2_PrThm2_Eq1}
	\sum_{i=[(k+5)/3] \atop i\ odd}^{k/2-1} d_i A^{(i)}
	=
	0,
\end{equation}
	where $c_j^{(i)}$ are rational numbers as in (ii),
	that is,  
	$c_i^{(i)} \neq c_{k-i}^{(i)}$ and 
	$c_j^{(i)} = c_{k-j}^{(i)}$ if $j$ is an odd integer with $[(k+5)/3] \leq j  < k/2$ and $j \neq i$.
However,
	it should hold that $c_h^{(h)} = c_{k-h}^{(h)}$ because of (\ref{2_PrThm2_Eq1}) and $d_h\neq0$,
	which is a contradiction.
Thus the relations obtained in (ii) are linearly independent.
Since the numbers of
the generators in (i) and $\zeta(odd, odd)$'s with weight $k$
equal $(k/2 - \mathrm{dim}\  M_k)$ and $(k-2)/2 $ respectively,
the number of the nontrivial equations equals $(\mathrm{dim}\ M_k -1)$.
\qed \\ \\

\begin{example}
We introduce examples of the relations obtained in Theorem \ref{1.Theorem2}(ii) in case of $k=18,20,22,24$.
(See \cite[Introduction]{GKZ} for the case of $k = 12, 16$.)
The case of $k=24$ is the first case where two relations appear,
	otherwise one relation.
Let $S_3$ be the symmetric group of degree $3$,
	and we put
\begin{equation*}
	s_j(k_1,k_2,k_3)
	=
	\sum\limits_{\sigma\in S_3} (-1)^{ k_{\sigma(1)} } \binom{ k_{\sigma(2)}+k_{\sigma(3)}-j-1 }{ k_{\sigma(2)}-1 }  
\end{equation*}
	for a nonnegative integer $j$ and positive integers $k_1,k_2,k_3$,
	where $\binom{m}{n}=0$ if $m<n$.
In order to give the examples,
	the explicit relations for double zeta values stated in Lemma \ref{2_Lem} are required,
	and we write them; 
\begin{multline} \label{3_Exm1_Eq1}
	0
	=
	\sum_{\sigma\in S_3}
	\Biggl\{
		\frac{2}{3} s_0(k_{\sigma(1)}, k_{\sigma(2)}, k_{\sigma(3)}) \zeta(k)
		+
		s_1(k_{\sigma(1)}, k_{\sigma(2)}, k_{\sigma(3)})  \zeta(1,k-1)
		\\
		+
		\sum_{j=2}^{k-2} s_j(k_{\sigma(1)}, k_{\sigma(2)}, k_{\sigma(3)})
		\Bigl[
			\zeta(j,k-j) 
			- 
			\Bigl( \delta_{j\ odd } + \frac{\delta_{j\ even}}{3} \Bigr) \zeta(j)\zeta(k-j)
		\Bigr]
	\Biggr\}
	,
\end{multline}	
	
\begin{eqnarray} \label{3_Exm1_Eq2}
		&&
		\left(\begin{array}{cccc}
		\binom{k-3}{1} & 0 & 0 & 0 \\
		\binom{k-3}{3} & \binom{k-5}{1} & 0 & 0 \\
		\binom{k-3}{5} & \binom{k-5}{3} & \binom{k-7}{1} & 0 \\
		\vdots &   &  & \ddots
		\end{array}\right)
		\left(\begin{array}{c}\zeta(2,k-2) \\\zeta(4,k-4) \\\zeta(6,k-6) \\\vdots\end{array}\right)
		\nonumber\\
		&=&
		-
		\left(\begin{array}{cccc}
		2 & 0 & 0 & 0 \\
		\binom{k-4}{2} & 2 & 0 & 0 \\
		\binom{k-4}{4} & \binom{k-6}{2} & 2 & 0 \\
		\vdots &   &  & \ddots
		\end{array}\right)
		\left(\begin{array}{c}\zeta(3,k-3) \\\zeta(5,k-5) \\\zeta(7,k-7) \\\vdots\end{array}\right)
		\nonumber\\
		&&
		-
		\left(\begin{array}{c}\binom{k-2}{2} \\\binom{k-2}{4} \\\binom{k-2}{6} \\\vdots\end{array}\right) \zeta(1,k-1)
		+
		\sum_{j=2}^{k-2} (-1)^{j-1} \left(\begin{array}{c}\binom{k-j-1}{2} \\\binom{k-j-1}{4} \\\binom{k-j-1}{6} \\\vdots\end{array}\right) \zeta(j)\zeta(k-j)
		,
\end{eqnarray}
	where 
	$k_1,k_2,k_3$ are positive integers with $k=k_1+k_2+k_3$,
	$\delta_{P(j)}$is the Kronecker delta function
	which equals $1$ if the condition $P(j)$ is true and $0$ if false,
	and
	$\zeta(1,k-1)$ means $-(\zeta(k) + \zeta(k-1,1))$.
Note that (\ref{3_Exm1_Eq1}) 
	can be derived from (\ref{2.Lemma1i}) and (\ref{2.Lemma1iib}) (\cite[Theorem 1.1.2]{nakamura1} and \cite[Corollary 2.4]{boyadzhiev1}),
	and
	(\ref{3_Exm1_Eq2}) from (\ref{2.Lemma1iii}) (\cite[Corollary 2.4]{boyadzhiev1} with $q=0$).
By the use of (\ref{3_Exm1_Eq1}) and (\ref{3_Exm1_Eq2}),
	we can calculate  the following examples.
\mbox{}\\
{\bf The case of $k=18$.}
By (\ref{3_Exm1_Eq1}) with $(k_1,k_2,k_3)=(6,6,6)$ and (\ref{3_Exm1_Eq2}), we obtain
\begin{multline*}
	4004\zeta(15,3)+23199\zeta(13,5)+47880\zeta(11,7)+59822\zeta(9,9)
	\\+
	47685\zeta(7,11)+24024\zeta(5,13)
	=
	\frac{62821831}{43867}\zeta(18)
	.
\end{multline*}
{\bf The case of $k=20$.}
By (\ref{3_Exm1_Eq1}) with $(k_1,k_2,k_3)=(8,6,6)$ and (\ref{3_Exm1_Eq2}), we obtain
\begin{multline*}
	858\zeta(17,3)+5005\zeta(15,5)+10758\zeta(13,7)+14925\zeta(11,9)
	\\
	+14938\zeta(9,11)+10725\zeta(7,13)+5148\zeta(5,15)
	=
	\frac{111230333}{349222}\zeta(20)
	.
\end{multline*}
{\bf The case of $k=22$.}
By (\ref{3_Exm1_Eq1}) with $(k_1,k_2,k_3)=(8,7,7)$ and (\ref{3_Exm1_Eq2}), we obtain
\begin{multline*}
	21216\zeta(19,3)+124566\zeta(17,5)+277732\zeta(15,7)+415239\zeta(13,9)
	\\
	+470415\zeta(11,11)+415324\zeta(9,13)+277290\zeta(7,15)+127296\zeta(5,17)
	\\
	=
	\frac{632571863}{77683}\zeta(22)
	.
\end{multline*}
{\bf The case of $k=24$.}
By (\ref{3_Exm1_Eq1}) with $(k_1,k_2,k_3)=(8,8,8), (10,7,7)$  and (\ref{3_Exm1_Eq2}), we obtain
\begin{multline*}
	1518100\zeta(21,3)+8953662\zeta(19,5)+20466050\zeta(17,7)+32180239\zeta(15,9)
	\\
	+39585975\zeta(13,11)+39585975\zeta(11,13)+32182500\zeta(9,15)
	\\
	+
	20447430\zeta(7,17)+9108600\zeta(5,19)
	=
	\frac{282358595588279}{472728182} \zeta(24)
	,
\end{multline*}
\begin{multline*}
	814606\zeta(21,3)+4807413\zeta(19,5)+11023896\zeta(17,7)+17421586\zeta(15,9)
	\\
	+21521181\zeta(13,11)+21522150\zeta(11,13)+17421586\zeta(9,15)
	\\+
	11015745\zeta(7,17)+4887636\zeta(5,19)
	=
	\frac{151914058887111}{472728182} \zeta(24)
	.
\end{multline*}
Note that,
	in the first equation,
	the coefficients of $\zeta(15,9)$ and $\zeta(9,15)$ are not equal,
	and the ones of $\zeta(13,11)$ and $\zeta(11,13)$ are equal.
On the other hand, in the second equation,
	the ones of $\zeta(15,9)$ and $\zeta(9,15)$ are equal,
	and the ones of $\zeta(13,11)$ and $\zeta(11,13)$ are not equal. 
\end{example}



\end{document}